\documentclass[11pt]{article}
\usepackage{amsmath,amsfonts,amssymb}
\usepackage{geometry}
\def\({\begin{eqnarray}}
\def\){\end{eqnarray}}
\def\[{\begin{eqnarray*}}
\def\]{\end{eqnarray*}}
\def\fref#1{(\ref{#1})}
\def\dt{\frac{d}{dt}}
\def\dx{\partial_x}
\def\intR{\int_{\R}}
\def\intRd{\int_{\R^d}}
\def\intRtwod{\int_{\R^d\times\R^d}}
\def\intRthreed{\int_{\R^d\times\R^d\times\R^d}}
\def\R{\mathbb{R}}
\def\N{\mathbb{N}}
\def\ind{{\rm 1\kern-.3em I}}

\newtheorem{Lem}{Lemma}
\newtheorem{The}{Theorem}
\newtheorem{Cor}{Corollary}

\def\begproof{\noindent{\bf Proof.}\ }
\def\endproof{\quad\vrule height4pt width4pt depth0pt \medskip}

\title{Aggregated Steady States of a Kinetic Model for Chemotaxis}

\author{Anne Nouri\footnote{Aix-Marseille University, CNRS, Centrale Marseille, I2M UMR 7373, 13453 Marseille, France,
{\tt anne.nouri@univ-amu.fr}} \ and Christian Schmeiser\footnote{University of Vienna,
Oskar-Morgenstern-Platz 1, 1090 Wien, Austria, {\tt Christian.Schmeiser@univie.ac.at}}
}

\begin{document}
\date{}\maketitle

\begin{abstract}
A kinetic chemotaxis model with attractive interaction by quasistationary chemical 
signalling is considered. The special choice of the turning operator, with velocity
jumps biased towards the chemical concentration gradient, permits closed ODE systems
for moments of the distribution function of arbitrary order. The system for second 
order moments exhibits a critical mass phenomeneon. The main result is existence
of an aggregated steady state for supercritical mass.
\end{abstract}

\section{Introduction}
Chemotaxis, the movement of biological agents influenced by gradients of chemical
concentrations, is a ubiquitous process in biological systems. On the other hand,
the production or degradation of chemicals is at the basis of standard signalling 
mechanisms between individuals. This produces a nonlinear feedback which, 
together with chemotactic motility, may drive self-organization processes in groups 
of agents. 

A typical example, observed in many bacterial and amoeboid species, 
is aggregation driven by the production of a diffusible chemical, and chemotactic
movement biased towards the direction of the gradient of the chemical concentration
(see the large literature on {\it Dictyostelium discoideum} or, for bacteria, \cite{Mittal}).
Since motility usually has a random component, it is a standard question in this 
situation, if the attractive mechanism is strong enough to overcome the dispersion
caused by the random motility component.

The type of mathematical models mostly depends on the nature of the 
motility process. The standard assumption of Brownian motion with drift, the latter
determined by chemotaxis, leads to a version of the classical Patlak-Keller-Segel
(PKS) model \cite{Patlak,KelSeg}, where a convection-diffusion equation for the agent 
density is coupled with a reaction-diffusion equation for the chemical concentration.
For certain bacterial species a description by a velocity jump process is more 
appropriate, whence the convection-diffusion equation of the PKS model is replaced
by a kinetic transport equation \cite{OthDunAlt}. The PKS model can typically
be recovered as a macroscopic limit \cite{CMPS,HwaKanSte}. However, some observed
phenomena are only explainable by kinetic models \cite{Saragosti}.

Three types of long time behavior can be observed in mathematical models. If the
random motion of agents dominates, this leads to dispersion, i.e. the same
qualitative behavior as for the heat equation. For dominating attractive effects, the
agent density either has a nontrivial aggregated long-time limit, or it blows up in finite 
time, typically in a concentration event. The two-dimensional parabolic-elliptic PKS
model (i.e. with a quasistationary equation for the chemical concentration) has been
thoroughly analyzed with respect to these questions. It shows a critical mass
phenomenon: Among the initial data with finite variance those with the total mass
below a critical value lead to dispersion and those with supercritical mass to finite
time blow-up \cite{BlaDolPer}. At the blow-up time strong solutions cease to exist,
but a continuation by measure solutions is possible as limiting case of regularized 
models \cite{DolSch,Vel}. The corresponding dichotomy has been shown to exist
also in kinetic transport models \cite{BouCal}. The situation for other versions of the
PKS model and, in particular, for kinetic models is less clear. 

Motivated by experimental results for {\it E. coli} \cite{Mittal}, a linear kinetic model
with given aggregated chemical concentration has been analyzed in \cite{CalRaoSch}.
The existence of a nontrivial steady state and its dynamic stability have been proven
(the latter by employing the methodology of \cite{DolMouSch}). The present work
can be seen as a continuation, where the nonlinear coupling with a quasistationary
model for the chemical is added. The main result is a critical mass phenomenon,
but with a dichotomy between dispersion and the existence of an aggregated steady
state.
Consider the system
\(
  && \partial_t f + v\dx f = \intR \bigl( T[S](v'\to v,x,t)f' - T[S](v\to v',x,t)f \bigr) dv' \,,
                                                       \label{kinetic}\\
  && -D\dx^2 S = \beta\rho_f - \gamma S \,,\label{elliptic}
\)
a one-dimensional kinetic model for chemotaxis, where the cells with phase space density $f(x,v,t)$
and macroscopic density and flux,
\[
  \rho_f(x,t) = \intR f(x,v,t) dv\quad\mbox{and, respectively,}\quad j_f(x,t) = \intR v f(x,v,t) dv\,,
\]
produce the chemoattractant with density $S(x,t)$. The dynamics of the chemoattractant (diffusion, production,
and decay) is assumed to be fast (and therefore modelled as quasistationary) compared to the dynamics of the
cells. 

We consider two choices for the turning kernel:
\[
  \mbox{\bf Model A:}\quad && T[S](v\to v',x,t) = \kappa S(x+\alpha v',t) \,,\\
  \mbox{\bf Model B:}\quad &&  T[S](v\to v',x,t) = \kappa S(x+\alpha(v'-v),t)\,.
\]
For both models, we assume $\alpha,\beta,\gamma,\kappa,D > 0$.
In Model A, cells decide about reorientation by scanning the chemoattractant density in the directions 
of possible post-turning velocities. In Model B, they scan in the direction of possible velocity changes. Both models have not been derived systematically from the 
microscopic behavior of a particular cell type. However, they are reasonable from a qualitative point of view, and they have the remarkable mathematical property that 
the evolution of moments can be computed by solving linear constant coefficient ODEs without solving the 
full equations (see Section \ref{sec:moments}).

We observe that the rescaling 
\[
  t\to \alpha t\,, \quad v\to \frac{v}{\alpha\sqrt{\gamma/D}}\,, \quad x \to \frac{x}{\sqrt{\gamma/D}}\,,\quad
  f\to f\,\frac{\alpha\gamma^2}{\kappa\beta D}\,, \quad 
  S\to S\,\frac{\sqrt{\gamma/D}}{\kappa}\,,
\] 
eliminates all parameters, i.e., \fref{kinetic}, \fref{elliptic} becomes
\(
  && \partial_t f + v \dx f = \intR \bigl( T[S](v'\to v,x,t)f' - T[S](v\to v',x,t)f \bigr) dv' \,,
                                                       \label{kinetic-s}\\
  && -\dx^2 S = \rho_f - S \,,\label{elliptic-s}
\)
with
\[
  \mbox{\bf Model A:}\quad && T[S](v\to v',x,t) = S(x+v',t) \,,\\
  \mbox{\bf Model B:}\quad &&  T[S](v\to v',x,t) = S(x+v'-v,t)\,.
\]

We consider the Cauchy problem with initial conditions
\( \label{IC}
  f(x,v,0) = f_I(x,v) \ge 0 \quad\mbox{for } x,v\in\R \,.
\)
The solution $S$ of \fref{elliptic} is defined as the convolution product of the decaying fundamental 
solution of $-\dx^2 + id$ with $\rho$: 
\(\label{equ:S}
   S[\rho](x,t) = \frac{1}{2} \intR e^{-|x-y|}\rho(y,t)dy \,.
\)
The initial datum is assumed to possess moments of up to second order:
\( 
  &&\intR \rho_I \,dx = \intR \intR f_I \,dv\,dx = M < \infty \,,\nonumber\\
  &&\intR |x|^2 \rho_I \,dx = \intR \intR |x|^2 f_I \,dv\,dx < \infty \,,\nonumber\\
  &&\intR\intR |v|^2 f_I \,dv\,dx < \infty \,.\label{Ass1}
\)
We choose a reference frame such that the first order moments vanish initially:
\(
  &&\intR j_I \,dx = \intR\intR vf_I \,dv\,dx = 0 \,,\nonumber\\
  && \intR x\rho_I \,dx = \intR\intR xf_I \,dv\,dx = 0 \,.\label{Ass2}
\)
We shall show that for Model A the ($x$- and $v$-) moments of $f$ up to any
fixed order satisfy closed systems of linear, constant coefficient ODEs. The system of second order moments
exhibits a critical mass phenomenon. If the total mass $M$ is below a critical value, the second order moments grow
indefinitely with time, whereas for large enough mass they converge to finite values. The corresponding system
for Model B always produces growing second order moments. Therefore we shall concentrate on model A
after this observation. It turns out that also the higher order moment systems exhibit a critical mass
phenomenon, however with the critical mass increasing with the moment order. Since stationary solutions may be the limits when time tends to infinity of the solutions to the Cauchy problem associated to (\ref{kinetic})-(\ref{elliptic}), this suggests a mass 
dependent decay behavior of the steady state. However, a precise characterization is still open.

The plan of the paper is the following. Global in time solutions to the Cauchy problem are determined in Section  \ref{cauchy-pb} for models A and B. The long term behavior of their moments is studied in Section  \ref{sec:moments} and proven to depend on the total mass of the solution.  A formal asymptotics of the solution to Model A is performed for large mass in Section \ref{formal-asympt}. The existence of a smooth steady state for Model A with supercritical mass (of the second order moment system) is proven in Section \ref{sec:stationary}.

\section{The Cauchy problem}\label{cauchy-pb}

\begin{The}\label{th1}
Given $f_I\in L^1_+(\R ^2)$ and $M= \int_{\R^2} f_I(x,v)dx\,dv$, there is a unique solution \\
$f\in \mathcal{C}([0,\infty),L^1_+(\R ^2))$ to the Cauchy problem associated to Model A, i.e.,
\begin{equation}\label{evolutionary-eqn}
\partial _tf+v\dx f = Q_A(f) \,,\qquad f(t=0)=f_I \,,
\end{equation}
with
\begin{equation}\label{evolutionary-eqn2}
Q_A(f)(v,x) = S[\rho_f](x+v)\rho_f(x) - M f(x,v) \,,\qquad  
S[\rho](x) = \frac{1}{2}\int_{\R} \rho (y)e^{-|x-y|}dy \,,
\end{equation}
where $\rho _f(x,t) = \int_{\R} f(x,v,t)dv$.
\end{The}

\begproof For every $T>0$, let
\begin{eqnarray*}
X_T := \left\{ \rho \in \mathcal{C}\left([0,T];L^1_+(\R )\right):\  \int_{\R} \rho(x,t)dx= M,\ \forall\, t\in [0,T]\right\} \,,
\end{eqnarray*}
equipped with the natural norm $\|\cdot\|_{X_T}$, and let $R_A(\rho )= \int_{\R} f\,dv$, where $f$ is the solution of
\begin{equation}\label{evolutionary-eqn3}
\partial _tf(x,v,t)+v\partial _xf(x,v,t)= \rho (x,t)S[\rho](x + v,t) - M f(x,v,t) \,,
\end{equation}
\begin{equation}\label{evolutionary-eqn4}
f(x,v,0)= f_I(x,v) \,.
\end{equation}
It can be computed explicitly as
\begin{equation}\label{R(rho)}
R_A(\rho)(x,t) = e^{-Mt}\int_{\R} f_I(x-vt,v)dv + \int_0^t e^{-Ms} \int_{\R} S[\rho](x+v(1-s),t-s)\rho(x-vs,t-s)dv\,ds.
\end{equation}
For $\rho\in X_T$ nonnegativity of $R_A(\rho)$ is obvious, and the mass conservation property follows by integration of \fref{R(rho)} with respect to $x$, implying $R:\, X_T\to X_T$.

The idea is to show that $R_A$ is a contraction with respect to $\|\cdot\|_{X_T}$. First we observe that
$\int_{\R} S[\rho](x,t)dx = M$ for $\rho\in X_T$ and that $\rho\mapsto S[\rho]$ as a map from $L^1(\R)$ 
to $L^1(\R)$ is Lipschitz with Lipschitz constant $1$. This implies for $\rho_1,\rho_2\in X_T$,
after a change of variables,
\begin{eqnarray*}
  \|R_A(\rho_1) - R_A(\rho_2)\|_{X_T} &\le&\int_0^T e^{-Ms} \int_{\R^2} \bigl(\rho_1(\xi,t-s)|S[\rho_1]-S[\rho_2]|
      (\eta,t-s)  \\
  && \qquad\qquad\qquad+ S[\rho_2](\eta,t-s) |\rho_1 - \rho_2|(\xi,t-s)\bigr)d\xi\,d\eta\,ds \\
   &\le& 2(1-e^{-MT})\|\rho_1 - \rho_2\|_{X_T} \,.
\end{eqnarray*}
Thus, for $T< \frac{\ln 2}{M}$ the map $R_A$ is a contraction on $X_T$. This proves local solvability. As a consequence of the uniform bound in $L^1(\R^2)$ the solution can be extended indefinitely in time steps 
of length $T$.
\endproof
\begin{The}\label{th1B}
Given $f_I\in L^1_+(\R ^2)$ and $M= \int_{\R^2} f_I(x,v)dx\,dv$, there is a unique solution \\
$f\in \mathcal{C}([0,\infty),L^1_+(\R ^2))$ to the Cauchy problem associated to Model B, i.e.,
\begin{equation}\label{evolutionary-eqnB}
\partial _tf+v\dx f = Q_B(f) \,,\qquad f(t=0)=f_I \,,
\end{equation}
with
\begin{equation}\label{evolutionary-eqn6}
Q_B(f)(v,x) = \int S[\rho_f](x+v-v^\prime )f(x,v^\prime )dv^\prime - M f(x,v) \,,\qquad  
S[\rho](x) = \frac{1}{2}\int_{\R} \rho (y)e^{-|x-y|}dy \,,
\end{equation}
where $\rho _f(x,t) = \int_{\R} f(x,v,t)dv$.
\end{The}

\begproof  A solution $f$ to the Cauchy problem associated to Model B is directly obtained as the limit of the increasing sequence $(f_j)$ defined by $f_0= 0$ and $f_{j+1}$ given from $f_j$ as the solution to
\begin{align}\label{fm-each-step}
&\partial _tf_{j+1}(x,v,t)+v\partial _xf_{j+1}(x,v,t)= \int S[\rho _{f_j}](x + v-v^\prime ,t)f_j(x,v^\prime,t)dv^\prime - M f_{j+1}(x,v,t) \,,\\
&f_{j+1}(x,v,0)= f_I(x,v) \, , \nonumber
\end{align}
where $\rho _{f_j}= \int f_jdv$. 
$f_{j+1}$ is explicitly given from $f_j$ by
\begin{equation}\label{fj}
f_{j+1}(x,v,t)= e^{-Mt}f_I(x-vt,v) + \int_0^t e^{M(s-t)} \int_{\R} S[\rho _{f_j}](x+v(s-t+1)-v^\prime ,s)f_j(x+v(s-t),v^\prime ,s)dv^\prime \,ds.
\end{equation}
Consequently it can be proven by induction that $(f_j)$ is nonnegative, and non decreasing since $f_0= 0 \leq f_1$ and $f_{j-1}\leq f_j$ imply $0\leq S[\rho _{f_{j-1}}] \leq S[\rho _{f_j}]$ and $f_j\leq f_{j+1}$ by \fref{fj}. Moreover, denoting by $m_j(t)= \int f_j(x,v,t)dxdv$ and integrating \fref{fm-each-step} with respect to $(x,v)\in \R ^2$, it holds 
\begin{eqnarray*}
m_{j+1}^\prime = m_j^2-Mm_{j+1}\, ,
\end{eqnarray*}
so that it can be proven by induction that 
\begin{equation}\label{bound-massB}
m_j(t)\leq M,\quad j\in \N .
\end{equation}
And so, by the monotone convergence theorem, $(f_j)$ converges in $L^1$ to a nonnegative function $f$. It implies that the limit in $L^1$ of $\big( \int S[\rho _{f_j}](x+v-v^\prime ,t)f_j(x,v^\prime ,t)dv^\prime \big) $ is $\int S[\rho _{f}](x+v-v^\prime ,t)f(x,v^\prime ,t)dv^\prime $. And so, $f$ is a solution of \fref{evolutionary-eqnB}. Conservation of mass from \fref{evolutionary-eqnB} implies that $\int f(x,v,t)dxdv= M$.
The limit of \fref{fj} implies that $f\in \mathcal{C}([0,\infty),L^1_+(\R ^2))$. $f$ is the unique nonnegative solution of \fref{evolutionary-eqnB} since its construction makes it minimal among the nonnegative solutions of \fref{evolutionary-eqnB}. Indeed, if there were another nonnegative solution $\tilde{f}$, then $f\leq \tilde{f}$ and $\int f(x,v,t)dxdv= \int \tilde{f}(x,v,t)dxdv= M$ would imply that $\tilde{f}= f$. \endproof

\section{Evolution of moments}\label{sec:moments}

The closedness of the equations for the moments for Model A relies on the following result

\begin{Lem}\label{lem-mom}
Let $0\le n\le N$, let $\rho\in L^1_+(\R)$ have finite moments up to order $N$, i.e., 
$$
  \intR |x|^k\rho(x)dx < \infty \,,\qquad k = 0,\ldots,N \,,
$$
and let $S$ be the bounded solution of $\dx^2 S = S-\rho$, i.e.,
$$
   S(x) = \frac{1}{2} \intR e^{-|x-y|}\rho(y)dy \,.
$$
Then also $S$ has finite moments up to order $N$, and with
$$
   R_k := \intR x^k \rho(x)dx\,,\qquad S_k := \intR x^k S(x)dx \,,\qquad k = 0,\ldots,N \,,
$$
the following relations hold:
\[
  \mbox{a)} && S_k = R_k + k(k-1)S_{k-2} \,,\qquad  k = 0,\ldots,N\,,\\
  \mbox{b)} && \intR \intR x^{N-n} v^n S(x+v)\rho(x)dx\, dv = \sum_{k=0}^n {n\choose k} (-1)^{n-k} S_k
                         R_{N-k} \,\quad 0\leq n\leq N.
\]
\end{Lem}

\begproof
The result can be shown by straightforward computations. We multiply the differential equation for $S$ by
$|x|^k$ and $x^k$, and use
$$
  \intR |x|^k \dx^2 S\,dx = k(k-1)\intR |x|^{k-2} S\,dx \,,\quad\mbox{and}\quad 
    \intR x^k \dx^2 S\,dx = k(k-1)S_{k-2} \,,
$$
to show the boundedness of the moments of $S$ and a). Then b) is a consequence of the substitution $y=x+v$
and of the binomial theorem.
\endproof

If we concentrate on moments of order $N$, then the lemma implies
\[
  S_N = R_N + LOT \,,
\]
and
\[
   \intR \intR x^{N-n} v^n S(x+v)\rho(x)dx\, dv = \left( (-1)^n  + \delta_{n,N} \right)R_0 R_N + LOT \,,
\]
where $LOT$ (Lower Order Terms) stands for terms only depending on moments of order lower than $N$.
Now we introduce moments of solutions $f$ of \fref{kinetic-s}, \fref{elliptic-s} with respect to $x$ and $v$:
$$
  A_{m,n}(t) := \intR\intR x^m v^n f(x,v,t)dv\,dx \,.
$$
With the help of Lemma \ref{lem-mom} and with 
$$
  -\intR\intR x^{N-n}v^n (v\dx f)dv\,dx = (N-n)A_{N-n-1,n+1} \,,
$$
we can derive differential equations for the moments. The first
and obvious one is mass conservation:
$$
  \dot A_{0,0} = 0 \qquad\Longrightarrow\qquad A_{0,0} = M \,.
$$
As a consequence, the turning operator of Model A after elimination of the unknown $S$ by \fref{equ:S}, can now be written as 
$$
  Q_A(f)(x,v) = S[\rho_f](x+v)\rho(x) - M f(x,v) \,.
$$
For the first order moments, we obtain
$$
  \dot A_{1,0} = A_{0,1} \,,\quad\dot A_{0,1} = -MA_{0,1} \qquad\Longrightarrow\qquad A_{1,0} = A_{0,1} = 0 \,,
$$
because of \fref{Ass2}. For the moments of order two it gets more interesting:
\[
  \dot A_{2,0} &=& 2 A_{1,1} \,,\\
  \dot A_{1,1} &=& A_{0,2} - M A_{1,1} - M A_{2,0} \,,\\
  \dot A_{0,2} &=& 2M A_{2,0} - M A_{0,2} + 2M^2\,.
\]
Application of the Routh-Hurwitz criterion to the characteristic polynomial of the coefficient matrix shows
that for $M<2$ at least one positive eigenvalue exists, whereas for $M>2$ all eigenvalues have negative real
parts. Thus, in the latter case all solutions converge to the steady state
$$
  (A_{2,0},A_{1,1},A_{0,2}) = \left( \frac{2M}{M-2}, 0, \frac{2M^2}{M-2}\right) \,.
$$
For the higher order moments we only concentrate on the highest order terms on the right hand sides:
\[
  \dot A_{N-n,n} = (N-n) A_{N-n-1,n+1} + \left( (-1)^n  + \delta_{n,N} \right)M A_{N,0} - M A_{N-n,n} + LOT\,,
\]
for $0\le n\le N$. This is a linear ODE system with constant coefficients and an inhomogeneity only depending
on lower order moments. This shows that all moments can be computed recursively.

If the coefficient matrices of all systems up to order $N$ only have eigenvalues with negative real parts,
then all moments of order up to $N$ have finite limits as $t\to\infty$. We have already shown above 
that for $N=2$ this property holds, iff $M>2$.

The characteristic polynomial of the order $N$ coefficient matrix can be written as
\[
  p_N(\lambda) = -\lambda (-M-\lambda)^N + M N!\sum_{n=0}^{N-1} \frac{(-M-\lambda)^n}{n!}
   + (-1)^N M N! \,,
\]
and the determinant of the coefficient matrix is thus given by
\[
  p_N(0) = (-1)^N M N! q_N(M) \,,\qquad q_N(M) = 1+ \sum_{n=0}^{N-1} (-1)^{N-n}\frac{M^n}{n!} \,.
\]
As we know already and can also be seen from $q_2(M) = 2-M$, the $2^{nd}$-order coefficient matrix has
a zero eigenvalue for $M=2$. The same is true for $3^{rd}$ order ($q_3(M) = M-M^2/2= (2-M)M/2$) but,
surprisingly, not for $4^{th}$ order. The function $q_4(M) = 2-M+M^2/2-M^3/6$ has a unique zero 
$M_4\in(2,3)$.\smallskip

\noindent{\bf Conjecture:} The functions $q_N$ have unique positive zeroes $M_N$, building an increasing
sequence, which tends to infinity. The $N^{th}$-order linear system above is stable, iff $M>M_N$.\smallskip

If the conjecture is true then, for every fixed $M>0$, only a finite number of moments tends to a bounded value
as $t\to\infty$. This would indicate an $M$-dependent decay of the equilibrium distribution with stronger
decay for larger values of $M$.

The essential parts of the conjecture can be proved:

\begin{Lem}
For fixed $N$, and $M$ large enough, all roots of $p_N$ have negative real parts.
\end{Lem}

\begproof
First we look for eigenvalues, which remain bounded as $M\to\infty$. For fixed $\lambda$,
\[
  \frac{p_N(\lambda)}{M^N} = (-1)^{N-1}(\lambda+N) + O(M^{-1}) \,,
\]
which provides a first root 
\[
  \lambda_0 = -N + O(M^{-1}) \,.
\]
Next we look for roots $\lambda=-M-\mu$ with $\mu$ bounded as $M\to\infty$. It is straightforward to show
\[
  \frac{p_N(-M-\mu)}{N! M} = r_N(\mu) + O(M^{-1}) \,,\qquad\mbox{with } 
    r_N(\mu) = (-1)^N + \sum_{n=0}^N \frac{\mu^n}{n!}  \,.
\]
Denoting the roots of $r_N$ by $\mu_1,\ldots,\mu_N\in\mathbb{C}$ (multiple roots allowed), we found 
$N$ more roots of $p_N$:
\[
  \lambda_j = -M-\mu_j + O(M^{-1}) \,,\qquad j=1,\ldots,N \,.
\]
Obviously, all the $N+1$ roots we found have negative real parts for large enough $M$.
\endproof

\begin{Lem}
For fixed $M$, and $N$ large enough, there exists at least one positive root of $p_N$.
\end{Lem}

\begproof
For fixed $\lambda$ and $M$,
\[
  \frac{p_N(\lambda)}{N! M} \approx e^{-M-\lambda} + (-1)^N \,,\qquad\mbox{as } N\to\infty \,.
\]
Therefore, for $N$ large enough there exists $\lambda>0$ such that sign$(p_N(\lambda)) = (-1)^N$.
On the other hand, 
\[
   \lim_{\lambda\to\infty} p_N(\lambda)(-1)^{N-1} = \infty \,, 
\]
completing the proof.
\endproof
\hspace*{0.1in}\\
Combination of the existence theorem \ref{th1} with the previous results leads to the propagation of moments:

\begin{Cor}
\hspace*{0.1in}\\
Let the assumptions of Theorem \ref{th1} hold and let $(1+|x|^N + |v|^N)f_I\in L^1(\R^2)$ for an $N\ge 1$. \\
Then the solution $f$ of \fref{evolutionary-eqn}, \fref{evolutionary-eqn2} satisfies 
$(1+|x|^N + |v|^N)f\in L^\infty_{loc}(\R^+;\,L^1(\R^2))$. \\
If $N\ge 2$ and $M>2$, then
 $(1+|x|^2 + |v|^2)f\in L^\infty(\R^+;\,L^1(\R^2))$.
\end{Cor}
\hspace*{0.1in}\\
For Model B, the computations are similar but a little more involved. As for Model A, $A_{0,0}=M$ and
$A_{1,0} = A_{0,1} = 0$ hold. The 2$^{nd}$ order moments satisfy the closed ODE system
\(
  \dot A_{2,0} &=& 2A_{1,1} \,,\nonumber\\
  \dot A_{1,1} &=& A_{0,2} - M A_{2,0} \,,\nonumber\\
  \dot A_{0,2} &=& 2M \left( M + A_{2,0} - A_{1,1} \right) \,.\label{ODEA}
\)
By their definition and by the Cauchy-Schwarz inequality, for nonvanishing initial data $f_I$ their initial values 
satisfy $A_{2,0}(0)A_{0,2}(0) > A_{1,1}(0)^2$ and $A_{2,0}(0), A_{0,2} > 0$. A straightforward computation gives
\[
  \dt (A_{2,0}A_{0,2}-A_{1,1}^2) = 2MA_{2,0}(M+A_{2,0}) \,.
\]
This guarantees that $A_{2,0}$ and $A_{0,2}$ remain positive for all times.

It is also easily seen that the Jacobian of the right hand side of \fref{ODEA} has at least one
positive eigenvalue. The only steady state has negative $A_{2,0}$- and $A_{0,2}$-components and can therefore
never be reached. Thus, all solution components tend to infinity exponentially,
meaning that the chemotactic effect is not strong enough to prevent dispersion.
For this reason we concentrate on Model A for the rest of this work.

\section{Formal asymptotics for large mass}\label{formal-asympt}

With the rescaling $f\to Mf$, $S\to MS$, Model A takes the form
\( \label{ModelA-macro}
  \partial_t f + v\dx f = M\left(S[\rho_f](x+v)\rho -  f\right) \,,
\)
with $M$ now taking the role of an inverse Knudsen number.

The rescaled version of the steady states for the moments are
\(\label{mom-inf}
    A_{2,0,\infty} = \frac{2}{M - 2} \,,\quad A_{1,1,\infty} = 0 \,,\quad 
  A_{0,2,\infty} = \frac{2M}{M - 2} \,,
\)
which suggests an equilibrium state concentrating with respect to $x$ as $M\to \infty$.

As $M\to \infty$, formally $f(x,v,t)\to f_0(x,v,t) = \rho_0(x,t)S[\rho_0](x+v,t)$. Mass conservation gives
\[
  \partial_t \rho_0 - \dx (x\rho_0) = 0 \,.
\]
Obviously, we have $\rho_0(x,t) \to \delta(x)$ as $t\to\infty$ and, thus,
$$
   \lim_{t\to\infty}f_0(x,v,t) = \frac{1}{2}e^{-|v|}\delta(x) \,. 
$$ 
This is in agreement with the limit as $M\to \infty$ in \fref{mom-inf}.
\section{Stationary solutions}\label{sec:stationary}
In this section we first prove in Theorem \ref{the:stationary} the existence of even nonnegative $L^1$solutions to the stationary problem, then their $C^\infty $ regularity in Theorem \ref{Cinfinity-regularity}.
\begin{The}\label{the:stationary}
For any $M>2$ there is an even nonnegative $L^1(\R ^2)$ solution $f$ of
\begin{equation}\label{stationary-model}
v\partial _xf(x,v)= \rho _f(x)S[\rho _f](x+v)-M f(x,v) \,,\qquad \int_{\R^2} f(x,v)dx\,dv= M \,.
\end{equation}
\end{The}
\begproof
As a consequence of Lemma \ref{lem-mom}, it holds that
\begin{equation}\label{S1}
\int_{\R} S[\rho](x+ v)dv= \int_{\R} \rho (y)dy \,,\qquad
  \int_{\R} vS[\rho](x+ v)dv= \int_{\R} \rho (y)(y-x)dy \,,
\end{equation}
\begin{equation}\label{S3}
\int_{\R} v^2S[\rho](x+v)dv= \int_{\R} \rho (y)((y-x)^2+2)dy \,. \hspace*{1.9in}
\end{equation}
Let $j\in \N^*$ and $M_j= M( 1-e^{-2j} - 1/j)$. Our first goal is to prove the existence and uniqueness of an even function $f_j\in L^1_+(\R^2)$, such that
\begin{equation}\label{eq1-fj}
f_j(x,v)= 0 \,,\qquad \mbox{for }|x| >j, \quad\text{or } |v| <\frac{1}{j}, \quad\text{or }|v| > 4j \,, \hspace*{1.3in}
\end{equation}
\begin{equation}\label{eq2-fj}
v\partial _xf_j(x,v)= \rho _{f_j}(x)S[\rho _{f_j}](x+v)- M f_j(x,v) \,,\qquad |x| <j \,,\quad \frac{1}{j}<|v| < 4j \,, 
\end{equation}
\begin{equation}\label{eq3-fj}
f_j(-j,v)= f_j(-j,-v),\quad \quad f_j(j,v)= f_j(j,-v), \hspace*{1.8in}
\end{equation}
\begin{equation}\label{eq4-fj}
\int_{\R^2} f_j(x,v)dxdv\in [ M_j, M] . \hspace*{3.2in}
\end{equation}
Let $K$ be the convex set 
\[ \begin{aligned}
K := \left\{ \rho \in L^1_+(\R):\ \rho \text{ even,}\ \rho (x)= 0\text{ if } |x| >j \,, \ \int_{\R} \rho (x)dx\in [ M_j, M] \right\} \,.
\end{aligned}\] 
Let the map $T$ be defined on $K$ by $T(\rho ) = \int_{\R} F\,dv$, where the restriction of $F$ to $\R^+\times \R $ is the solution of
\begin{equation}\label{eq1-F}
F(x,v)= 0 \,,\qquad\mbox{for } x>j \,,\quad\text{or } |v| <\frac{1}{j}\,, \quad\text{or } |v| >4j \,, \hspace*{1.4in}
\end{equation}
\begin{equation}\label{eq2-F}
v\partial _xF(x,v)= \rho (x)S[\rho](x+ v)-M F(x,v) \,,\qquad 0<x<j \,,\quad \frac{1}{j}< |v| <4j \,, \hspace*{0.1in}
\end{equation}
\begin{equation}\label{eq3-F}
F(0,v)= e^{Mj/v}\big( F(j,-v)-\int_{-j/v}^0e^{M\tau }\rho (j+\tau v)S[\rho](j+\tau v+ v)d\tau \big) \,,\quad v>0 \,,
\end{equation}
and $F$ is extended by parity w.r.t. $(x,v)$ to $\R^-\times \R $.
Equations (\ref{eq2-F}), (\ref{eq3-F}) imply that
\begin{equation}\label{eq4-F}
F(j,v)= F(j,-v).
\end{equation}
Denote by $X_j$ the space of nonnegative $L^1\left( \left (- 4j, -\frac{1}{j}\right) \right) $ functions with the weight $|v|$. \\
The solution $F\in L^1_+(\R ^+\times \R )$ of (\ref{eq1-F})--(\ref{eq3-F}) exists and is unique, because the map 
\begin{eqnarray*}
\gamma \in X_j\rightarrow F(j,-v) \,,\hspace*{0.05in} -j<v<-\frac{1}{j} \,,
\end{eqnarray*}
where $F$ is the solution of (\ref{eq1-F}), (\ref{eq2-F}) and 
\begin{equation}\label{eq5-F}
F(j,v)= \gamma (v),\hspace*{0.05in} v<0,\qquad F(0,v)= e^{Mj/v}\left( \gamma (-v)-\int _{-j/v}^0e^{M\tau }\rho (j+\tau v)S[\rho|(j+\tau v+v)d\tau \right) ,\hspace*{0.05in} v>0, \nonumber
\end{equation}
is a contraction. Indeed, for any $(\gamma _1,\gamma _2)\in X_j^2$ with images $(F_1(j,-v),F_2(j,-v))_{-4j<v<-\frac{1}{j}}$, 
we have
\begin{align*}
\int _{-j}^{-1/j} |v| \,|F_1(j,-v)-F_2(j,-v)| dv&= \int_{-4j}^{-1/j} |v| \, |\gamma_1(v)-\gamma _2(v)| e^{2Mj/v}dv\\ 
&\leq e^{-M/2}\int_{-j}^{-1/j} |v| \,|\gamma_1(v)-\gamma_2(v)| dv \,. 
\end{align*}
Due to the exponential form of \fref{eq2-F}, $F$ is nonnegative. Hence, $T(\rho)$ is nonnegative. Moreover, $T(\rho)$ is even since $F$ is even. Equation (\ref{eq2-F}) holds on $(-j,j) \times \big( (-4j,-1/j) \cup (1/j,4j) \big) $ since $F$, $\rho $ and $S_\rho $ are even functions. Integrating it on $(-j,j) \times \big( (-4j,-1/j) \cup (1/j,4j) \big) $ and using (\ref{eq1-F}) and (\ref{S1}), implies that 
\begin{eqnarray*}
M \int_{\R^2} F\,dx\,dv= M^2 - \int_{\R} \rho (x)\int _{|v| \in (0,1/j) \cup  (4j,\infty)}S[\rho ](x+ v)dv\,dx \,.
\end{eqnarray*}
Moreover,
\begin{eqnarray*}
\int_{|v| \in (0,1/j) \cup  (4j,\infty)}S[\rho ](x+ v)dv\leq \big(e^{-2j} + 1/j\big) \int _{-j}^j\rho (y)dy \,,\qquad |x| <j \,.
\end{eqnarray*}
And so, $T$ maps $K$ into $K$. We claim that
$T$ is compact with respect to the $L^1$ topology. Indeed, let $(\rho _n)$ be a sequence in $K$, i.e., such that 
\begin{eqnarray*}
\int_{\R} \rho _n(x)dx\in [ M_j, M] \,,\qquad n\in \N \,.
\end{eqnarray*}
By definition of $K$, the sequence $(T(\rho_n))= \left(\int_{\R} F_ndv\right)$ satisfies
\begin{eqnarray*}
\int_{\R} T(\rho _n)(x)dx\in [ M_j,M] \,.
\end{eqnarray*}
Moreover,
\begin{eqnarray*}
 |\partial _xF_n(x,v)| <jM\big( \rho _n(x)+F_n(x,v)\big) \,,\qquad |x| <j \,,\quad \frac{1}{j}< |v| < 4j \,,
\end{eqnarray*}
so that 
\begin{eqnarray*}
\int_{-j}^j |\partial_x T(\rho_n)(x)|dx \leq 2M^2j \,.
\end{eqnarray*}
The Sobolev space $W^{1,1}( (-j,j) )$ being compactly imbedded in $L^1( (-j,j) )$, there is a subsequence $(T(\rho_{n_k}))$ of $\left(T(\rho_n)|_{(-j,j)}\right)$ converging in $L^1( (-j,j) )$. Consequently the map $T$ is compact.

$T$ is continuous with respect to the $L^1$ topology. Indeed, let $(\rho _n)$ be a sequence in $K$ that converges to $\rho $ in $L^1(\R )$. By the previous compactness argument, there is a subsequence $(T(\rho_{n_k}))$ of 
$(T(\rho _n))$ converging in $L^1$ to some $\sigma = \int_{\R} F(x,v)dv$. $F$ solves (\ref{eq2-F}) since 
$(\rho _n)$ converges to $\rho $ in $L^1(\R )$ and $(S[\rho _n])$ converges to $S[\rho]$ in $C([ -j,j] )$. Moreover,  the solution $F$ of (\ref{eq2-F}), (\ref{eq3-F}) is unique. Indeed, if there were two solutions, their difference 
$G$ would satisfy
\begin{equation}\label{eq1-G}
v\partial _xG(x,v)= -M G(x,v) \,,\qquad |x| <j \,,\quad \frac{1}{j}< |v| < 4j \,,
\end{equation}
\begin{equation}\label{eq2-G}
G(\pm j,v)= G(\pm j,-v) \,.
\end{equation}
Multiplying (\ref{eq1-G}) by $G$, integrating the resulting equation over $[ -j,j] \times \{ 1/j< |v| < j\} $ and using (\ref{eq2-G}) implies that $\int_{\R^2} G^2(x,v)dxdv= 0$, i.e. $G$ is identically zero. And so, the whole sequence $(F_{n})$ converges in $L^1$ to $F$.
Thus, there is a fixed point $f_j$ of $T$, i.e. an even solution of (\ref{eq1-fj})--(\ref{eq4-fj}).\\
Prove that there is some constant $c$ such that
\begin{eqnarray*}
\int (1+x^2+v^2)f_j(x,v)dxdv<c,\quad j\in \N ^*.
\end{eqnarray*}
Multiplying (\ref{eq2-fj}) by $xv$  (resp. $v^2$, resp. $x^2$), integrating the corresponding equation on \\
$[-j,j] \times [-4j,-\frac{1}{j}] \cup [\frac{1}{j},4j] $ and using Lemma \ref{lem-mom} leads to
\begin{align*}
&\int_{\R^2} v^2f_j\,dx\,dv-M\int_{\R} x^2f_j\,dxdv= M\int xvf_jdxdv+j\int_{\R} v^2(f(j,v)+f(-j,v))dv\nonumber \\
&\hspace*{2.2in}+\int_{\R} x\rho _{f_j}(x)\int _{|v| \in [ 0,1/j ] \cup  [ 4j,\infty) }vS[\rho _{f_j}](x+ v)dv\,dx \,, \\
&\int v^2f_jdxdv-2\int x^2f_jdxdv= 2M-\frac{1}{M}\int_{\R} \rho _{f_j}(x) \int_{|v| \in [ 0,1/j] \cup  [4j,\infty)}v^2S[\rho _{f_j}](x+ v)dv\,dx \,, \\
&2\int xvf_jdxdv= \int_{\R} \rho _{f_j}(x) x^2\int_{|v| \in [ 0,1/j] \cup  [4j,\infty)}S[\rho _{f_j}](x+ v)dv\,dx \,.
\end{align*}
Hence,
\begin{align*}
(M-2)\int v^2f_j(x,v)ddv&= 2M^2-2\int_{\R} \rho _{f_j}(x)x \int_{|v| \in [ 0,1/j] \cup  [4j,\infty)}vS[\rho _{f_j}](x+ v)dv\,dx \\
&-2j\int v^2(f(j,v)+f(-j,v))dv-M\int_{\R} \rho _{f_j}(x)x^2 \int_{|v| \in [ 0,1/j] \cup  [4j,\infty)}S[\rho _{f_j}](x+ v)dv\,dx \\
&-\int_{\R} \rho _{f_j}(x)\int_{|v| \in [ 0,1/j] \cup  [4j,\infty)}v^2S[\rho _{f_j}](x+ v)dv\,dx \\
&\leq 2M^2-2\int_{-j}^j \rho _{f_j}(x)x \int_{|v| \in [ 0,1/j] \cup  [4j,\infty)}vS[\rho _{f_j}](x+ v)dv\,dx \, 
\end{align*}
and
\begin{align*}
(M-2)\int x^2f_j(x,v)ddv&= 2M-\int_{\R} \rho _{f_j}(x)x \int_{|v| \in [ 0,1/j] \cup  [4j,\infty)}vS[\rho _{f_j}](x+ v)dv\,dx \\
&-\frac{M}{2}\int_{\R} \rho _{f_j}(x)x^2 \int_{|v| \in [ 0,1/j] \cup  [4j,\infty)}S[\rho _{f_j}](x+ v)dv\,dx \\
&-\frac{1}{M}\int_{\R} \rho _{f_j}(x)\int_{|v| \in [ 0,1/j] \cup  [4j,\infty)}v^2S[\rho _{f_j}](x+ v)dv\,dx \\
&\leq 2M-\int_{-j}^j \rho _{f_j}(x)x \int_{|v| \in [ 0,1/j] \cup  [4j,\infty)}vS[\rho _{f_j}](x+ v)dv\,dx \, .
\end{align*}
Moreover, for $M>2$,
\begin{align*}
\lvert \int_{-j}^j \rho _{f_j}(x)x \int_{|v| \in [ 0,1/j] \cup  [4j,\infty)}vS[\rho _{f_j}](x+ v)dv\,dx \rvert &\leq \big( \frac{1}{j}+12je^{-2j}\big) M\int _{-j}^j \rho _{f_j}(x)\lvert x\rvert dx\\
&\leq 4M\int \lvert x\rvert f_j(x,v)dxdv\\
&\leq \frac{M-2}{2}\int x^2f_j(x,v)dxdv+\frac{8M^3}{M-2}\, .
\end{align*}
Finally, 
\begin{align}\label{moment-order-2}
\int_{\R^2} (1+x^2+v^2)f_j \,dx\,dv &\leq \frac{60M^3}{(M-2)^2} \,,\quad j\in \N ^* \, .
\end{align}
It follows from \fref{moment-order-2} that a subsequence of $(f_j)$ tightly converges to a nonnegative bounded measure $\mu $. We still denote this subsequence by $(f_j)$. The measure $\mu $ has its total mass equal to $M$ like every $f_j$. In order to prove that $\mu $ satisfies the first equation of \fref{stationary-model}, we start from \fref{eq2-fj} written along the characteristics, multiply it by a continuous and bounded test function $\varphi (x,v)$, so that
\begin{align}\label{fj-weak-form}
\int f_j(x,v)\varphi (x,v)dxdv= &\int _{v>0}\int _{\R }\int _{-\infty }^0 \rho _{f_j}(x+sv)S_{\rho _{f_j}}(x+(1+s)v)e^{Ms}\varphi (x,v)dsdxdv\nonumber \\ 
+&\int _{v<0}\int _{\R }\int _0^{+\infty }\rho _{f_j}(x-sv)S_{\rho _{f_j}}(x+(1-s)v)e^{-Ms}\varphi (x,v)dsdxdv\nonumber \\
= &\int \rho _{f_j}(x)\Big( \int _{v>0}\int _{-\infty }^0 S_{\rho _{f_j}}(x+v)e^{Ms}\varphi (x-sv,v)dsdv\Big) dx\nonumber \\
+&\int \rho _{f_j}(x)\Big( \int _{v<0}\int _0^{+\infty } S_{\rho _{f_j}}(x+v)e^{-Ms}\varphi (x+sv,v)dsdv\Big) dx.
\end{align}
The left hand side of \fref{fj-weak-form} tends to $<\mu ,\varphi >$ when $j\rightarrow +\infty $ because of the tight convergence of $(f_j)$ to $\mu $. The right hand side of \fref{fj-weak-form} tends to 
\begin{align*}
&\int \rho _{\mu }(x)\Big( \int _{v>0}\int _{-\infty }^0 S_{\rho _{\mu }}(x+v)e^{Ms}\varphi (x-sv,v)dsdv\Big) dx\\
+&\int \rho _{\mu }(x)\Big( \int _{v<0}\int _0^{+\infty } S_{\rho _{\mu }}(x+v)e^{-Ms}\varphi (x+sv,v)dsdv\Big) dx\\
&=\int _{v>0}\int _{\R }\int _{-\infty }^0 \rho _{\mu }(x+sv)S_{\rho _{\mu }}(x+(1+s)v)e^{Ms}\varphi (x,v)dsdxdv\\
&+\int _{v<0}\int _{\R }\int _0^{+\infty }\rho _{\mu }(x-sv)S_{\rho _{\mu }}(x+(1-s)v)e^{-Ms}\varphi (x,v)dsdxdv
\end{align*}
 when $j\rightarrow +\infty $, because $(\rho _{f_j})$ tends to $\rho _{\mu }$ in $L^\infty (\R )$ weak star and 
 \begin{eqnarray*} 
 x\rightarrow \int _{v>0}\int _{-\infty }^0S_{\rho _{f_j}}(x+v)e^{Ms}\varphi (x-sv,v)dsdv\quad \Big( \text{resp.   } x\rightarrow \int _{v<0}\int _0^{+\infty }S_{\rho _{f_j}}(x+v)e^{-Ms}\varphi (x+sv,v)dsdv\Big)
\end{eqnarray*}  
converges in $L^1(\R )$ to
\begin{eqnarray*} 
 x\rightarrow \int _{v>0}\int _{-\infty }^0S_{\rho _{\mu }}(x+v)e^{-Ms}\varphi (x-sv,v)dsdv\quad \Big( \text{resp.   } x\rightarrow \int _{v<0}\int _0^{+\infty }S_{\rho _{\mu }}(x+v)e^{-Ms}\varphi (x-sv,v)dsdv\Big)
\end{eqnarray*}
by the dominated convergence theorem. And so, $\mu $ is a nonnegative bounded measure stationary solution to the problem.\\
Moreover, integrating (\ref{eq3-fj}) between $-\infty $ and $x$ (resp. $x$ and $+\infty $) for $v>0$ (resp. $v<0$) implies that
\begin{eqnarray*}
 |v| f_j(x,v)\leq M\int_{\R }\rho _{f_j} dx = M^2\,,\qquad (x,v)\in \R ^2,\quad j\in \N \,.
\end{eqnarray*}
Consequently the only singular part of the measure $\mu $ is a Dirac measure at $v= 0$. Let us split $\mu $ as
\begin{eqnarray*}
\mu = g+\gamma \delta(v) \,,
\end{eqnarray*}
with $g\in L^1(\R ^2)$ and $\gamma \in L^1(\R )$. Equation (\ref{eq3-fj}) for $\mu $ writes
\begin{eqnarray*}
0= (\rho_g + \gamma )S[\rho_g +\gamma](x+v) - Mg - M\gamma \,\delta(v) \,, 
\end{eqnarray*}
so that $\gamma = 0$. It follows that $\mu \in L^1(\R^2)$.
\endproof

\begin{The}\label{Cinfinity-regularity}
Let $M>2$ hold. Then solutions $f$ of \eqref{stationary-model} as in Theorem \ref{the:stationary} satisfy
$f\in C^\infty(\R^2)$.
\end{The}

\begproof
The solution of the approximative problem \eqref{eq1-fj}--\eqref{eq4-fj} in the proof of Theorem
\ref{the:stationary} satisfies
\[
  f_j(x,v) = f_j(-j,v)e^{-M(x+j)/v} + \int_0^{(x+j)/v} \rho_{f_j}(x-sv)S[\rho_{f_j}](x+v(1-s)) e^{-Ms}ds \,,
\]
for $v>0$, $-j\le x\le j$, and
\[
  f_j(x,v) = f_j(j,v)e^{-M(x-j)/v} + \int_0^{(x-j)/v} \rho_{f_j}(x-sv)S[\rho_{f_j}](x+v(1-s)) e^{-Ms}ds \,,
\]
for $v<0$, $-j\le x\le j$. The estimates in the proof allow to pass to the limit $j\to\infty$,
showing that the problem \eqref{stationary-model} is solved in the mild sense:
\(\label{mild-stat}
  f(x,v) = \int_0^{+\infty }\rho_f (x-sv)S[\rho_f](x+v(1-s)) e^{-Ms}ds \,,\quad (x,v)\in \R ^2.
\)
With the Fourier transform
\[
  \hat f(\xi,k) = \intR \intR f(x,v) e^{-i(\xi\cdot x + k\cdot v)} dv\,dx \,,
\]
a straightforward computation leads to
\begin{eqnarray*}
\widehat{S[\rho]}(\xi )= \frac{\widehat\rho (\xi )}{1+\xi ^2}\,,\quad \xi \in \R .
\end{eqnarray*}
Consequently,  \eqref{mild-stat} is equivalent to
\(\label{mild-stat-F}
  \widehat{f}(\xi,k) = \int_0^{+\infty }e^{-Ms}
    \frac{\widehat{\rho_f}(\xi(1+s)-k)\widehat{\rho_f}(k-\xi s)}{1+(k-\xi s)^2} ds \,.
\)
Moreover, $\widehat{\rho_f}(\xi) =   \widehat f(\xi,0)$, so that
\[
  \widehat{\rho_f}(\xi)= \int_0^{+\infty }e^{-Ms}
  \frac{\widehat{\rho_f}(\xi(1+s))\widehat{\rho_f}(-\xi s)}{1+\xi^2 s^2} ds \,.
\]
Using the boundedness of $\widehat{\rho_f}$ by $M$ on the right hand side leads to
$\widehat{\rho_f}(\xi) = O(|\xi|^{-1})$ as $|\xi|\to\infty$. This can be iterated,
giving $\widehat{\rho_f}(\xi) = O(|\xi|^{-n})$ for arbitrary $n$ and, therefore,
$\rho_f\in C^\infty(\R)$. Actually, we shall use 
\[
  |\widehat{\rho_f}(\xi)| \le \frac{c_n}{(1+\xi^2)^n} \,,\qquad \forall n\ge 0 \,,
\]
in \eqref{mild-stat-F}, leading to the estimate
\[
 |\widehat{f}(\xi,k)| \le c_n^2 \int_0^{+\infty } 
 \frac{e^{-Ms} }{(1+(\xi(1+s)-k)^2)^n (1+(\xi s-k)^2)^n} ds \,.
\]
Note that
\[
  (1+(\xi(1+s)-k)^2) (1+(\xi s-k)^2) \ge 1+\max\{(\xi(1+s)-k)^2, (\xi s-k)^2\} \,.
\]
Minimizing the right hand side with respect to $k$ and, respectively, with respect 
to $\xi$, we obtain
\[
  (1+(\xi(1+s)-k)^2) (1+(\xi s-k)^2) &\ge& 1 
   + \max\left\{ \frac{\xi^2}{4}, \frac{k^2}{(2s+1)^2}\right\} \\
  &\ge& 1 + \frac{\xi^2}{5} + \frac{k^2}{5(2s+1)^2} \,.
\]
This shows that 
\[
  \frac{1 + \xi^2 + k^2}{(2s+1)^2 (1+(\xi(1+s)-k)^2) (1+(\xi s-k)^2)} \le
  \frac{1 + \xi^2 + k^2}{(2s+1)^2 (1+\xi^2/5) + k^2/5} \le 5 \,,
\]
and, consequentially, for every $n\ge 0$ there exists $C_n>0$, such that
\[
  |\widehat{f}(\xi,k)| \le \frac{5^n c_n^2}{(1+\xi^2+k^2)^n} 
  \int_0^\infty e^{-Ms} (2s+1)^{2n} ds = \frac{C_n}{(1+\xi^2+k^2)^n} \,,
\]
implying $f\in C^\infty(\R^2)$.
\endproof

\end{document}